\documentclass[12pt]{amsart}

\usepackage{amssymb,amsmath}
\usepackage[all]{xy}
\usepackage{setspace}

\usepackage{amsmath}
\usepackage{amssymb}
\usepackage{url}
\usepackage{amscd}
\usepackage{mathrsfs}
\usepackage[dvips]{graphicx}

\usepackage[OT2,T1]{fontenc}
\DeclareSymbolFont{cyrletters}{OT2}{wncyr}{m}{n}
\DeclareMathSymbol{\Sha}{\mathalpha}{cyrletters}{"58}

\theoremstyle{plain}
\newtheorem{thm}{Theorem}

\newtheorem{lemma}{Lemma}
\newtheorem{prop}{Proposition}

\newcommand{\Z}{{\mathbb Z}}
\newcommand{\N}{{\mathbb N}}

\newcommand{\R}{{\mathbb R}}

\newcommand{\CC}{{\mathbb C}}
\newcommand{\F}{{\mathbb F}_q}
\newcommand{\E}{{\mathbb E}}

\newcommand{\jac}{\mathrm{Jac}}

\newcommand{\J}{\mathcal J}

\newcommand{\HH}{{\mathcal H}}

\newcommand{\ud}{\mathrm{d}}

\begin{document}


\title[Jacobians of hyperelliptic curves over finite fields]{Statistics of the Jacobians of hyperelliptic curves over finite fields}



\author[Xiong]{Maosheng Xiong}
\author[Zaharescu]{Alexandru Zaharescu$^{\star}$}

\address{Maosheng Xiong: Department of Mathematics \\
Hong Kong University of Science and Technology \\
Clear Water Bay, Kowloon\\
P. R. China}
\email{mamsxiong@ust.hk}

\address{Alexandru Zaharescu: Institute of Mathematics of the Romanian Academy \\
P.O. Box 1--764, 70700 Bucharest, Romania \\
and Department of Mathematics \\
University of Illinois at Urbana-Champaign \\
273 Altgeld Hall, MC-382 \\
1409 W. Green Street \\
Urbana, IL 61801 USA}
\email{zaharesc@math.uiuc.edu}



\thanks{$^{\star}$The second author was supported by NSF grant number DMS-0901621.}

\begin{abstract}
Let $C$ be a smooth projective curve of genus $g \ge 1$ over a finite field $\F$ of cardinality $q$. In this paper, we first study $\#\J_C$, the size of the Jacobian of $C$ over $\F$ in case that $\F(C)/\F(X)$ is a geometric Galois extension. This improves results of Shparlinski \cite{shp}. Then we study fluctuations of the quantity $\log \#\J_C-g \log q$ as the curve $C$ varies over a large family of hyperelliptic curves of genus $g$. For fixed genus and growing $q$, Katz and Sarnak showed that $\sqrt{q}\left(\log \# \J_C-g \log q\right)$ is distributed as the trace of a random $2g \times 2g$ unitary symplectic matrix. When the finite field is fixed and the genus grows, we find the limiting distribution of $\log \#\J_C-g \log q$ in terms of the characteristic function. When both the genus and the finite field grow, we find that $\sqrt{q}\left(\log \# \J_C-g \log q\right)$ has a standard Gaussian distribution.

\end{abstract}

\maketitle

\section{Introduction}
Let $C$ be a smooth projective curve of genus $g \ge 1$ over a finite field $\F$ of cardinality $q$. The Jacobian $\jac(C)$ is a $g$-dimensional abelian variety. The set of the $\F$-rational points on $\jac(C)$, denoted by $\J_C=\jac(C)(\F)$, is a finite abelian group. The group $\J_C$ has been studied extensively, partly because of its importance in the theory of algebraic curves and its surprising applications in public-key cryptography and computational number theory. For example, such groups are extremely useful in primality testing \cite{adl} and integer factorization \cite{len1,len2}. Statistics of group structures of $\J_C$, for instance the analogue of the Cohen-Lenstra conjecture over function fields remains an inspiring problem in number theory and provides insight for number fields case. Interested readers may refer to \cite{ach1,ach2,ven} for details and current development. The main purpose of this paper is to study $\#J_C$, the size of the Jacobian over $\F$. This quantity is also the class number of the function field $\F(C)$ (\cite[Theorem 5.9]{ros}), a subject of study with a rich history.

The zeta function of $C/\F$ is a rational function of the form
\[ Z_{C}(u) = \frac{P_{C}(u)}{
(1-u)(1-qu)}\,, \]
where $P_{C}(u) \in \Z[u]$ is a polynomial of degree $2g$ with $P_{C}(0)=1$, satisfying the functional equation
\[P_{C}(u)=\left(qu^2\right)^g P_{C}\left(\frac{1}{qu}\right), \]
and having all its zeros on the circle $|u| = 1/\sqrt{q}$ (the Riemann Hypothesis for curves \cite{wei}). Moreover, there is a unitary symplectic matrix $\Theta_C \in \mathrm{USp}(2g)$, defined up to conjugacy, so that
\[P_C(u) =\det \left(I-u \sqrt{q} \Theta_C\right)\,.\]
The eigenvalues of $\Theta_C$ are of the form $e\left(\theta_{C,j}\right), j=1, \ldots, 2g$, where $e(\theta)=e^{2 \pi i \theta}$.

It is known that $\# \J_C=P_C(1)$ (see \cite[Corollary VIII.6.3]{lor}). From this we immediately derive that
\[(q^{1/2}-1)^{2g} \le \# \J_C \le (q^{1/2}+1)^{2g},\]
which is tight in the case $g=1$ due to the classical result of Deuring \cite{deu}. Many improvements of this bound have been obtained in \cite{que,rosen,shp,ste,tsf}. In particular in an interesting paper \cite{shp}, Shparlinski proves that if $C$ is a smooth absolutely irreducible curve of genus $g$ over $\F$ with gonality $d$, then
\begin{eqnarray} \label{1:shp}
\log \# \J_C= g \log q+O\left(g \log^{-1}(g/d)\right)
\end{eqnarray}
as $g \to \infty$, where the implied constant may depend on $q$. (The gonality of a curve $C$ is the smallest integer $d$ such that $C$ admits a non-constant map of degree $d$ to the projective line over the ground field $\F$. For example, a hyperelliptic curve is a curve given by an affine model $Y^2=F(X)$ for some $F \in \F[X]$, so the gonality is $d = 2$.) This generalizes and improves similar results of Tsfasman \cite{tsf}.

In this paper, we first prove that in case the function field $\F(C)$ is a geometric Galois extension of $\F(X)$, a sharper estimate can be obtained. Here ``geometric'' means that the constant field of $\F(C)$ is still $\F$. 

\begin{thm} \label{1:thm1} Let $C$ be a smooth projective curve of genus $g \ge 1$ over $\F$. Assume that the function field $\F(C)$ is a geometric Galois extension of the rational function field $\F(X)$ with $N=\#\mathrm{Gal}\left(\F(C)/F[X]\right)$. Then
\begin{eqnarray} \label{1:ine}
|\log \# \J_C- g \log q| \le \left(N-1\right)\left( \log \max \left\{1, \frac{\log(7g/(N-1))}{\log q}\right\}+3\right)\,.
\end{eqnarray}
\end{thm}

We remark that in Theorem 1, the inequality (\ref{1:ine}) is explicit and holds true for any $g$ and $q$. Moreover, the quantity $\log \# \J_C- g \log q$ is essentially bounded by $O(\log \log g)$, which is significantly smaller than $O(g/\log g)$ implied from (\ref{1:shp}).

Now assume that $q$ is odd. For each positive integer $d \ge 3$, denote by $\HH_{d,q}$ the family of hyperelliptic curves having an affine equation of the form $Y^2=F(X)$, with $F \in \F[X]$ a monic square-free polynomials of degree $d$. The measure on $\HH_{d,q}$ is simply the uniform probability measure on the set of such polynomials. The genus of a curve $C \in \HH_{d,q}$ is given by
\[g=g(C)=\left[\frac{d-1}{2}\right],\]
where for $x \in \R$, $[x]$ denotes the largest integer not exceeding $x$. For any $C \in \HH_{d,q}$, since $\F(C)/\F(X)$ is a geometric Galois extension with Galois group $\Z/2\Z$, Theorem \ref{1:thm1} implies that
\begin{eqnarray} \label{1:hy} \left|\log \# \J_{C}-g \log q\right| \le \log \max \left\{1,\log \frac{\log (7g)}{\log q}\right\}+3\,.\end{eqnarray}
The inequality (\ref{1:hy}) appears to be very sharp, however we will see that for most hyperelliptic curves $C \in \HH_{d,q}$, the value $|\log \# \J_{C}-g \log q|$ is actually much smaller. More precisely we prove the following. 

\begin{thm} \label{1:thm2}
For any $\psi \ge 2$, denote
\[M_{\psi}:=\#\left\{C \in \HH_{d,q}: |\log \#\J_{C}-g \log q| \ge \psi\right\}. \]
Then
\begin{eqnarray*} \frac{\#M_{\psi}}{\# \HH_{d,q}} \ll \exp\left(-\frac{\psi}{2} \log(q\log \psi)\right),\end{eqnarray*}
where the implied constant in ``$\ll$'' is absolute.
\end{thm}

Theorem \ref{1:thm2} shows that the ratio $\frac{\#M_{\psi}}{\# \HH_d}$ goes to zero very fast whenever $\psi$ or $q$ goes to infinity. When $q \to \infty$, this is not surprising. Actually, in this case, much more is known. Writing \[P_C(u)=\prod_{i=1}^{2g}\left(1-\sqrt{q} e\left(\theta_{C,i}\right)u\right),\]
then
\begin{eqnarray*} \log \#\J_{C}-g \log q=\sum_{i=1}^{2g} \log \left(1-q^{-1/2}e(\theta_{C,i})\right). \end{eqnarray*}
Katz and Sarnak \cite{kat} showed that for fixed genus $g$, the conjugacy classes $\left\{\Theta_{C}: C\in \HH_{d,q}\right\}$ become uniformly distributed in $\mathrm{USp}(2g)$ in the limit $q \to \infty$. In particular, since
\[\lim_{q \to \infty}\sqrt{q}\left(\log \#\J_{C}-g \log q\right)=-\sum_{i=1}^{2g} e(\theta_{C,i}),\]
it implies that
\begin{itemize}
\item[(i).] When $g$ is fixed and $q \to \infty$, the value $-\sqrt{q}\left(\log \#\J_{C}-g \log q\right)$ for $C \in \HH_{d,q}$ is distributed asymptotically as the trace of a random matrix in $\mathrm{USp}(2g)$.
\end{itemize}
Furthermore, since the limiting distribution of traces of a random matrix in $\mathrm{USp}(2g)$, as $g \to \infty$, is a standard Gaussian by a theorem of Diaconis and Shahshahani \cite{dia}, it also implies that
\begin{itemize}
\item[(ii).] If $q \to \infty$ and then $g \to \infty$, the value $\sqrt{q}\left(\log \#\J_{C}-g \log q\right)$ is distributed as a standard Gaussian.
\end{itemize}
Katz and Sarnak's powerful theorem \cite{kat} provides an almost complete story, except that in their argument, it is crucial to take the limit that $q \to \infty$. What happens if $g \to \infty$ instead? Complementary to (i) and (ii) above, we prove the following.

\begin{thm} \label{1:thm3} (1). If $q$ is fixed and $g \to \infty$, then for $C \in \HH_{d,q}$, the quantity \\
$\log \#\J_{C}-g \log q+\delta_{d/2}\log \left(1-q^{-1}\right)$ converges weakly to a random variable $X$, whose characteristic function $\phi(t)=\E\left(e^{itX}\right)$ is given by
\[\phi(t)= 1+\sum_{r=1}^{\infty} \frac{1}{r!} \sum_{\substack{P_1,\ldots,P_r\\ \mbox{\tiny distinct}}} \prod_{j=1}^r \frac{\left(1-|P_j|^{-1}\right)^{-it}+
\left(1+|P_j|^{-1}\right)^{-it}-2}{2\left(1+|P_j|^{-1}\right)}, \quad \forall t \in \R,\]
where
\[\delta_{\gamma}=\left\{\begin{array}{cc}1&\gamma \in \Z,\\
0& \gamma \not \in \Z, \end{array}\right.\]
and the sum on the right over $P_1,\ldots,P_r$ is over all distinct monic irreducible polynomials $P_1,\ldots,P_r \in \F[X]$ and $|P_j|=q^{\deg P_j}$.

(2). If both $q,g \to \infty$, then for $C \in \HH_{d,q}$, $\sqrt{q}\left(\log \#\J_{C}-g \log q\right)$ is distributed as a standard Gaussian, that is, for any $\gamma \in \R$, we have
\[\lim_{\substack{q \to \infty\\ g \to \infty}} \frac{1}{\# \HH_{d,q}} \#\left\{C \in \HH_{d,q}: \sqrt{q}\left(\log \#\J_{C}-g \log q\right) \le \gamma\right\}=\frac{1}{\sqrt{2 \pi}} \int_{-\infty}^{\gamma} e^{-\frac{t^2}{2}}\, \ud t\,.\]
\end{thm}

We remark that first, Theorem \ref{1:thm3} is in the spirit of Kurlberg and Rudnick \cite{kur} and Faifman and Rudnick \cite{fai}, who initiated the investigation of such problems. It is also in the spirit of Bucur, David, Feigon and Lal\'in \cite{buc1,buc2}, who made further important development. The proof of Theorem \ref{1:thm3} depends heavily on techniques developed by Rudnick \cite{rud} and Faifman and Rudnick \cite{fai}. Second, (2) of Theorem \ref{1:thm3} is slightly more general than Statement (ii) above as there is no requirement that $q \to \infty$ first. Finally, instead of averaging over $\HH_{d,q}$, the family of hyperelliptic curves arising from monic square-free polynomials of degree $d$, the proof of Theorem \ref{1:thm3} can be easily adapted to the moduli space of hyperelliptic curves of a fixed genus. Interested readers may refer to \cite{buc1,buc2} for terminology and treatment.

We organize this paper as follows. In Section 2 we collect several results which will be used later. Then we prove Theorem 1 in Section 3, Theorem 2 in Section 4 and Theorem 3 in Section 5 respectively. The proofs of Theorem 2 and Theorem 3 rely on Proposition 1--3 which are of technical natural. To emphasis and streamline the main ideas, we deal with these three propositions in Section 6.


\section{Preliminaries}
In this section we collect several results which will be used later. Interested readers can refer to \cite{ros} for more details.

\subsection{Zeta functions of function fields}

Let $K=\F(X)$ be the rational function field over the finite field $\F$ and let $L/K$ be a finite geometric Galois extension. Here ``geometric'' means that the constant field of $L$ is still $\F$. We list several facts about such extensions $L/K$ as follows (see \cite[Chapter 9]{ros} for more details).

First, the zeta function $\zeta_L(s)$ of $L$ is defined by
\[\zeta_L(s)=\prod_{P \in \mathcal{S}_L} \left(1-|P|^{-s}\right)^{-1},\]
where the product is over $\mathcal{S}_L$, the set of all primes of $L$, and for each $P \in \mathcal{S}_L$, $|P|$ is the cardinality of the residue field of $L$ at $P$. For the rational function field $K$, the zeta function $\zeta_K(s)$ turns out to be
\begin{eqnarray*}  \zeta_K(s)= \left(1-q^{-s}\right)^{-1} \left(1-q^{1-s}\right)^{-1}\,.\end{eqnarray*}
If $C$ is a smooth projective curve of genus $g \ge 1$ over $\F$ with function field $\F(C)=L$, then $Z_{C} \left(q^{-s}\right)=\zeta_L(s)$, i.e., the zeta function of the curve $C$ coincides with the zeta function of the function field $\F(C)$ (see \cite[pp. 57, Chap 5]{ros} for details).

Let $G=\mathrm{Gal}(L/K)$ be the Galois group of $L/K$ and $\rho: G \to \mathrm{Aut}_{\CC}(V)$ a representation of $G$, where $V$ is a finite-dimensional vector space over the complex numbers $\CC$ of dimension $m$. One defines the Artin L-series associated to the representation $\rho$ as follows.

If $P$ is a prime of $K$ which is unramified in $L$ and $\mathcal{B}$ is a prime of $L$ lying above $P$, one defines the local factor $L_P(s,\rho)$ as
\begin{eqnarray} \label{2:lf}L_P(s,\rho)=\det \left(I-\rho((\mathcal{B},L/K))|P|^{-s}\right)^{-1}\,,\end{eqnarray}
where $I$ is the identity automorphism on $V$ and $(\mathcal{B},L/K) \in G$ is the Frobenius automorphism at $\mathcal{B}$. Since $L/K$ is Galois, this definition does not depend on the choice of $\mathcal{B}$ over $P$.

Let $\{\alpha_1(P), \alpha_2(P), \ldots, \alpha_m(P)\}$ be the eigenvalues of $\rho((\mathcal{B},L/K))$. In terms of these eigenvalues, we get another useful expression for $L_P(s,\rho)$:
\[L_P(s,\rho)^{-1}=\left(1-\alpha_1(P)|P|^{-s}\right)\left(1-\alpha_2(P)|P|^{-s}\right) \cdots \left(1-\alpha_m(P)|P|^{-s}\right)\,.\]
We note that these eigenvalues $\alpha_i(P)$ are all roots of unity because $(\mathcal{B},L/K)$ has finite order.

At a prime $P$ of $K$ which is ramified in $L$, the local factor $L_P(s,\rho)$ can also be defined. The definition is similar to (\ref{2:lf}), except that the action $\rho((\mathcal{B},L/K))$ is restricted to a subspace of $V$ which is fixed by the inertial group $I(\mathcal{B}/P)$. We are contended with the fact that there are only finitely many primes $P$ which are ramified in $L$ and in either case we can write $L_P(s,\rho)$ as
\[L_P(s,\rho)^{-1}=\left(1-\alpha_1(P)|P|^{-s}\right)\left(1-\alpha_2(P)|P|^{-s}\right) \cdots \left(1-\alpha_m(P)|P|^{-s}\right)\,,\]
where the values $\alpha_i(P)$'s are either roots of unity or zero. The Artin L-series $L(s,\rho)$ is defined by the infinite product
\[L(s,\rho)=\prod_{P \in \mathcal{S}_K} L_P(s,\rho), \]
where $\mathcal{S}_K$ is the set of all primes in $K=\F(X)$.

It is known that if $\rho=\rho_0$, the trivial representation, then $L(s,\rho_0)=\zeta_K(s)$, and if $\rho=\rho_{\mathrm{reg}}$, the regular representation, then $L(s,\rho_{\mathrm{reg}})=\zeta_L(s)$. It is also known that $L(s,\rho)$ depends only on the character $\chi$ of $\rho$, so we can write it as $L(s,\chi)$. 

Finally, let $L/K$ be a finite, geometric and Galois extension with Galois group $G=\mathrm{Gal}(L/K)$. Let $\{\chi_1,\chi_2,\ldots,\chi_h\}$ be the set of irreducible characters of $G$. We set $\chi_1=\chi_0$, the trivial character. Denote by $d_i$ the degree of $\chi_i$, i.e., $d_i=\chi_i(e)$ is the dimension of the representation space corresponding to $\chi_i$. Then using results about group characters and formal properties of Artin L-series, one derives that
\begin{eqnarray} \label{2:lfunction}
\zeta_L(s)=\zeta_K(s) \prod_{i=2}^h L(s,\chi_i)^{d_i}\,.
\end{eqnarray}

\subsection{Averaging over $\HH_{d,q}$} 

Let $\HH_{d,q} \subset \F[X]$ be the set of all monic square-free polynomials of degree $d \ge 3$. 

\begin{lemma} \label{2:lem1} For any Dirichlet character $\chi: \F[X] \to \CC$ modulo $f \in \F[X]$, we have
\[\frac{1}{\#\HH_{d,q}} \sum_{F \in \HH_{d,q}} \chi(F) \le \frac{2^{\deg f-1}}{\left(1-q^{-1}\right)q^{d/2}}\,.\]
\end{lemma}
\noindent {\bf Proof.} This is \cite[Lemma 3.1]{fai}, which proves the case when $\chi=\left(\frac{f}{\cdot}\right)$ is a quadratic character. For the general case, the proof follows exactly the same line of argument, so we omit the details here. \quad $\square$

\begin{lemma} \label{2:lem2} Let $h \in \F[X]$ be a monic square-free polynomial. Then
\[\frac{1}{\#\HH_{d,q}} \sum_{\substack{F \in \HH_{d,q}\\ \gcd(F,h)=1}} 1 = \prod_{P|h}\left(1+|P|^{-1}\right)^{-1}+O\left(q^{-d/2} \sigma(h)\right)\,,\]
where $\sigma(h)=\sum_{D|h}1$.
\end{lemma}
\noindent {\bf Proof.} This is essentially \cite[Lemma 5]{rud}, which treats the case that $h=P$ is a monic irreducible polynomial. In fact in this case \cite[Lemma 5]{rud} yields a much stronger error term $O(q^{-d})$. The extra saving is obtained by carefully analyzing the functional equation of the zeta function. To get the error term $O(q^{-d/2}\sigma(h))$, the proof follows a standard procedure which is included \cite[Lemma 5]{rud}. We also omit details here. \quad $\square$

\section{Proof of Theorem \ref{1:thm1}}

Let $C$ be a smooth projective curve of genus $g \ge 1$ over $\F$. The zeta function $Z_C(u)$ is of the form
\[Z_C(u) =\frac{P_C(u)}{(1-u)(1-qu)}\,, \]
where $P_{C}(u) \in \Z[u]$ is a polynomial of degree $2g$ with $P_{C}(0)=1$, satisfying the functional equation
\[P_{C}(u)=\left(qu^2\right)^g P_{C}\left(\frac{1}{qu}\right), \]
and having all its zeros on the circle $|u| = 1/\sqrt{q}$. We may write $P_C(u)$ as
\[P_C(u)= \prod_{i=1}^{2g} \left(1-\sqrt{q}e(\theta_i) u\right),\]
where these $\theta_i \in [0,1)$ and $e(\alpha)$ stands for $e^{2 \pi i \alpha}$ for any $\alpha \in \R$.

Since $\#\J_C=P_C(1)$, we have
\begin{eqnarray*}  \# \J_C=\prod_{i=1}^{2g} \left(1-\sqrt{q}e(\theta_i)\right) =q^g \prod_{i=1}^{2g} \left(1-q^{-1/2}e(\theta_i) \right)\,.\end{eqnarray*}
Taking logarithms on both sides and using the expansion
\begin{eqnarray} \label{3:exp}-\log (1-z)=\sum_{n \ge 1} \frac{z^n}{n}, \quad |z|<1,\end{eqnarray}
we obtain the equation
\begin{eqnarray} \label{3:log}
\log \#\J_C-g \log q =\sum_{n \ge 1} q^{-n/2}n^{-1} \sum_{i=1}^{2g} -e(n \theta_i)\,.
\end{eqnarray}

Denote $L=\F(C)$ and $K=\F(X)$. The zeta functions of $L$ and $K$ can be written as
\[\zeta_L(s) =(1-q^{-s})^{-1}(1-q^{1-s})^{-1} \, \prod_{i=1}^{2g} \left(1-\sqrt{q}e(\theta_i) q^{-s}\right), \]
and
\[\zeta_K(s)=(1-q^{-s})^{-1}(1-q^{1-s})^{-1}\,.\]
Since $L/K$ is a geometric Galois extension with $G=\mathrm{Gal}(L/K)$ and $\# G=N$, let $\{\chi_1,\chi_2,\ldots,\chi_h\}$ be the set of irreducible characters of $G$ with $\chi_1=\chi_0$, the trivial character and denote by $d_i$ the degree of $\chi_i$. From (\ref{2:lfunction}) we find that
\begin{eqnarray} \label{3:lfunction} \prod_{i=2}^h L(s,\chi_i)^{d_i}=\prod_{i=1}^{2g} \left(1-\sqrt{q}e(\theta_i) q^{-s}\right),\end{eqnarray}
where for each $i$ with $2 \le i \le h$, the Artin L-series associated to $\chi_i$ can be written as
\[L(s,\chi_i)^{-1}=\prod_{P} \left(1-\alpha_{i,1}(P)|P|^{-s}\right)\left(1-\alpha_{i,2}(P)|P|^{-s}\right) \cdots \left(1-\alpha_{i,d_i}(P)|P|^{-s}\right)\,. \]
Here the product is over all monic irreducible polynomials $P \in \F(X)$ and $P=\infty$ with $|P|=q^{\deg P}$ ($\deg \infty=1$ hence $|\infty|=q$) and these $\alpha_{i,j}(P)$'s are either roots of unity or zero.

Taking logarithms on both sides of (\ref{3:lfunction}), using the expansion (\ref{3:exp}) again and equating the coefficients, we obtain for any positive integer $n$ the identity
\begin{eqnarray} \label{3:id} q^{n/2} \sum_{j=1}^{2g}-e(n \theta_i)=\sum_{\deg f=n} \Lambda(f) \sum_{i=2}^h d_i\sum_{j=1}^{d_i}\alpha_{i,j}(f)\,, \end{eqnarray}
where the sum on the right side over $\deg f=n$ is over all monic polynomials $f \in \F[X]$ with $\deg f=n$, $\Lambda(f)=\deg P$ if $f=P^k$ is a prime power, and $\Lambda(f)=0$ otherwise.

Let $Z$ be a positive integer which will be chosen later. Denote
\[\epsilon_{1,Z}=\sum_{n \le Z} q^{-n/2}n^{-1} \sum_{i=1}^{2g} -e(n \theta_i),\]
and
\[\epsilon_{2,Z}=\sum_{n > Z} q^{-n/2}n^{-1} \sum_{i=1}^{2g} -e(n \theta_i).\]
From (\ref{3:log}) we can write
\[\log \#\J_C-g \log q=\epsilon_{1,Z}+\epsilon_{2,Z}\,.\]
If $Z \ge 2$ we have
\begin{eqnarray} \label{3:e22} |\epsilon_{2,Z}| \le \sum_{n \ge Z+1}q^{-n/2}n^{-1}2g \le \frac{2g}{Z+1}\,q^{-(Z+1)/2}\left(1-q^{-1/2}\right)^{-1},\end{eqnarray}
and if $Z=1$ we have
\begin{eqnarray} \label{3:e21} |\epsilon_{2,Z}| \le 2g  \left(-\log\left(1-q^{-1/2}\right)-q^{-1/2}\right) \le \frac{2g}{q-\sqrt{q}}\,.\end{eqnarray}

For $\epsilon_{1,Z}$, we use the identity (\ref{3:id}). Since $|\alpha_{i,j}| \le 1$ for all $i,j$, we obtain the inequality
\[|\epsilon_{1,Z}| \le \sum_{n \le Z} q^{-n}n^{-1} \sum_{\deg f=n}\Lambda(f) \sum_{i=2}^h d_i^2\,.\]
It is known that
\[1+\sum_{i=2}^h d_i^2=N=\# G\]
and
\[\sum_{\deg f=n}\Lambda(f)=q^n+1.\]
Here the extra ``1'' on the right side in the above equation accounts for $f=\infty^n$. Hence
\[|\epsilon_{1,Z}| \le (N-1) \left(\sum_{n \le Z}\frac{1}{n}+\sum_{n \le Z}\frac{1}{nq^n}\right)\,.\]
If $Z=1$, this is
\begin{eqnarray} \label{3:e11} |\epsilon_{1,Z}| \le (N-1)\left(1+q^{-1}\right), \end{eqnarray}
and if $Z \ge 2$, we use
\[\sum_{n \le Z}\frac{1}{n} \le 1.5+\log Z-\log 2\]
and
\[\sum_{n \le Z}\frac{1}{nq^n} \le -\log \left(1-q^{-1}\right) \le \frac{1}{q-1}\]
to obtain
\begin{eqnarray} \label{3:e12} |\epsilon_{1,Z}| \le (N-1) \left(1.5-\log 2+\frac{1}{q-1}+\log Z\right),\quad Z \ge 2\,.\end{eqnarray}

Case 1: if $2\left(1-q^{-1/2}\right)^{-1}g \ge (N-1)q$, we choose
\[Z=\left[\frac{2 \log \frac{2\left(1-q^{-1/2}\right)^{-1}g}{N-1}}{\log q}\right] \ge 2\,.\]
We find from (\ref{3:e12}) that
\[|\epsilon_{1,Z}| \le (N-1) \left\{1.5+\frac{1}{q-1}+\log \left(\frac{\log \frac{2\left(1-q^{-1/2}\right)^{-1}g}{N-1}}{\log q}\right)\right\}\]
and from (\ref{3:e22}) that
\[|\epsilon_{2,Z}| \le \frac{N-1}{2}\,.\]
In this case noticing that $q \ge 2$, we obtain
\[\left|\log \#\J_C-g \log q\right|\le (N-1) \left(\log \left(\frac{\log \frac{7g}{N-1}}{\log q}\right)+3\right)\,.\]

Case 2: if $2\left(1-q^{-1/2}\right)^{-1}g < (N-1)q$, we choose $Z=1$, and from (\ref{3:e11}) and (\ref{3:e21}) we obtain that
\[\left|\log \#\J_C-g \log q\right|\le (N-1) \left(2+q^{-1}\right) < 3(N-1)\,.\]
In either case we conclude that
\[|\log \# \J_C- g \log q| \le \left(N-1\right)\left( \log \max \left\{1, \frac{\log(7g/(N-1))}{\log q}\right\}+3\right)\,.\]
This completes the proof of Theorem \ref{1:thm1}. \quad $\square$

\section{Proof of Theorem \ref{1:thm2}}

\subsection{Preparation} Let $\F$ be a finite field of cardinality $q$ with $q$ odd. Denote
\[\HH_{d,q}=\left\{F \in \F[X]: F \mbox{ is monic, square-free and} \deg F=d\right\}\,.\]
For any $F \in \HH_{d,q}$, the hyperelliptic curve $C_F$ is given by the affine model
\[C_F:Y^2=F(X).\]
It has genus
\[g=g_{_F}=\left[\frac{d-1}{2}\right]\,.\]
Suppose that the zeta function $Z_{C_F}(u)$ is of the form
\[Z_{C_F}(u) =\frac{\prod_{i=1}^{2g} \left(1-\sqrt{q}e\left(\theta_{i,F}\right) u\right)}{(1-u)(1-qu)}\,, \]
where the $\theta_{i,F}$'s are real numbers. Then
\begin{eqnarray*} \# \J_{C_F}=\prod_{i=1}^{2g} \left(1-\sqrt{q}e\left(\theta_{i,F}\right)\right) =q^g \prod_{i=1}^{2g} \left(1-q^{-1/2}e\left(\theta_{i,F}\right) \right)\,.\end{eqnarray*}
Taking logarithms on both sides we obtain the equation
\begin{eqnarray*}
\log \#\J_{C_F}-g \log q =\sum_{n \ge 1} q^{-n/2}n^{-1} \sum_{i=1}^{2g} -e\left(n \theta_{i,F}\right)\,.
\end{eqnarray*}
Let $Z$ be a positive integer which will be chosen later. 
We write
\begin{eqnarray} \label{4:jc}\log \#\J_{C_F}-g \log q=\sum_{n \le Z} q^{-n/2}n^{-1} \sum_{i=1}^{2g} -e\left(n \theta_{i,F}\right)+\epsilon_{1,Z}(F)\,,\end{eqnarray}
where
\[
\epsilon_{1,Z}(F)=\sum_{n > Z} q^{-n/2}n^{-1} \sum_{i=1}^{2g} -e\left(n \theta_{i,F}\right)\,.
\]
It is easy to see that
\begin{eqnarray*}
|\epsilon_{1,Z}(F)| \le \sum_{n > Z} q^{-n/2}n^{-1} 2g \le \frac{9g}{Z}\,\,q^{-Z/2}\,.
\end{eqnarray*}

Denote $L=\F(C_F)$ and $K=\F(X)$. Since $L/K$ is a geometric quadratic extension and the Legendre symbol $ \chi:=\left(\frac{F}{\cdot}\right)$ generates the Galois group Gal$(L/K)$, from (\ref{2:lfunction}) we have \begin{eqnarray} \label{4:lfunction} L\left(s,\chi\right)=\prod_{i=1}^{2g} \left(1-\sqrt{q}e\left(\theta_{i,F}\right) q^{-s}\right),\end{eqnarray}
and by definition
\begin{eqnarray} \label{4:lfunction2} L(s,\chi)=\prod_{P} \left(1-\left(\frac{F}{P}\right)|P|^{-s}\right)^{-1}\,. \end{eqnarray}
Here the product is over all monic irreducible polynomials $P \in \F(X)$ and $P=\infty$ with $|P|=q^{\deg P}$ ($\deg \infty=1$ hence $|\infty|=q$).

Computing $\frac{\ud}{\ud s} L(s,\chi)$ in two different ways using (\ref{4:lfunction}) and (\ref{4:lfunction2}) and equating the coefficients we obtain for each positive integer $n$ the identity
\begin{eqnarray} \label{4:fun} \sum_{i=1}^{2g}-e\left(n \theta_{i,F}\right)=q^{-n/2}\sum_{\deg f=n}\Lambda(f) \left(\frac{F}{f}\right)+q^{-n/2} \delta_{d/2}\,,\end{eqnarray}
where the sum over $\deg f=n$ on the right side is over all monic polynomials $f \in \F[X]$ with $\deg f=n$, and for any $\gamma \in \R$, $\delta_{\gamma}=1$ if $\gamma \in \Z$, and $\delta_{\gamma}=0$ if $\gamma \not \in \Z$. The extra term $q^{-n/2}\delta_{n/2}$ comes from $f=\infty^n$, noting the fact that $F \in \HH_{d,q}$ is monic and
\[\left(\frac{F}{\infty}\right)=\left\{\begin{array}{ccc}1&:& \deg F \equiv 0 \pmod{2}, \\
0&:& \deg F \equiv 1 \pmod{2}. \end{array}\right.\]
Using the identity (\ref{4:fun}) in (\ref{4:jc}) and denoting
\begin{eqnarray*} N_F=\log \# \J_{C_F}-g \log q+\delta_{d/2} \log \left(1-q^{-1}\right),\end{eqnarray*}
we find that
\begin{eqnarray*} N_F=\triangle_Z(F)+\epsilon_Z(F),\end{eqnarray*}
where
\begin{eqnarray} \label{4:tri}
\triangle_Z(F)=\sum_{n \le Z}q^{-n} n^{-1} \sum_{\deg f=n} \Lambda(f) \left(\frac{F}{f}\right),
\end{eqnarray}
and
\begin{eqnarray} \label{4:e}
|\epsilon_Z(F)| \le \frac{10g}{Z} \,\, q^{-Z/2}\,.
\end{eqnarray}

An upper bound for $\triangle_Z(F)$ is given by
\[\left|\triangle_Z(F)\right| \le \sum_{n \le Z}q^{-n}n^{-1} \sum_{\deg f=n} \Lambda(f)\le 1+\log Z\,,\]
and this can be used to estimate the quantity $N_F$ (as in the proof of Theorem \ref{1:thm1}). Extra savings can be obtained by taking high moments of $N_F$ and then averaging over the set $\HH_{d,q}$.

\subsection{The $r$-th moment $\triangle_Z$}

For any function $\chi: \HH_d \to \CC$, we denote by $\langle \chi \rangle$ the mean value of $\chi$ on $\HH_{d,q}$, that is,
\[\langle \chi \rangle:=\frac{1}{\# \HH_{d,q}} \sum_{F \in \HH_{d,q}} \chi(F)\,. \]
Assume that $d \ge 100$. We choose $r$ to be any positive even integer in the range
\begin{eqnarray} \label{4:r} 4 \le r \le \log d \,.\end{eqnarray}
We will estimate the $r$-moment $\left\langle (\triangle_Z)^r\right\rangle$ first.

From (\ref{4:tri}),
\[\triangle_Z(F)^r=\sum_{n_1,\ldots,n_r \le Z} \prod_{i=1}^r q^{-n_i}n_i^{-1} \sum_{\substack{\deg f_i=n_i\\ 1 \le i \le r}} \Lambda(f_1) \cdots \Lambda(f_r) \left(\frac{F}{f_1 \cdots f_r}\right)\,,\]
hence
\[\left\langle(\triangle_Z)^r \right\rangle=\sum_{n_1,\ldots,n_r \le Z} \prod_{i=1}^r q^{-n_i}n_i^{-1} \sum_{\substack{\deg f_i=n_i\\ 1 \le i \le r}} \Lambda(f_1) \cdots \Lambda(f_r) \left\langle\left(\frac{\cdot}{f_1 \cdots f_r}\right)\right \rangle \,.\]
If $f_1 \cdots f_r$ is not a square in $\F[X]$, then $\left(\frac{\cdot}{f_1 \cdots f_r}\right):\F[X] \to \CC$ is a non-trivial Dirichlet character modulo $h$ with $\deg h \le \sum_{i=1}^r \deg f_i$, by Lemma \ref{2:lem1} we find that
\[ \left\langle\left(\frac{\cdot}{f_1 \cdots f_r}\right)\right \rangle \le \frac{2^{n_1+\cdots+n_r-1}}{\left(1-q^{-1}\right)q^{d/2}}\,.\]
The total contribution to $\langle(\triangle_Z)^r \rangle$ from this case is bounded by
\[T_1 \le \sum_{n_1, \ldots, n_r \le Z} \prod_{i=1}^r q^{-n_i}n_i^{-1} \sum_{\substack{\deg f_i =n_i \\ 1 \le i \le r}} \Lambda(f_1) \cdots \Lambda(f_r) \frac{2^{n_1+\cdots+n_r-1}}{\left(1-q^{-1}\right)q^{d/2}}.\]
This can be estimated as
\begin{eqnarray} \label{4:t1} T_1 \le \frac{q^{-d/2} 2^{(Z+1)r}}{2\left(1-q^{-1}\right)} \le q^{-d/2} 2^{(Z+1)r}\,\, \,.\end{eqnarray}

If $f_1 \cdots f_r$ is a square in $\F[X]$, denote $f_1 \cdots f_r=h^2$ and $\widetilde{h}=\prod_{P|h}P$, then $\left(\frac{\cdot}{h^2}\right)$ is a trivial character, by Lemma \ref{2:lem2} we find that
\[ \left\langle\left(\frac{\cdot}{h^2}\right)\right \rangle = \frac{1}{\#\HH_{d,q}} \sum_{\substack{F \in \HH_{d,q}\\ \gcd(F,\widetilde{h})=1}}1= \prod_{P|\widetilde{h}}\left(1+|P|^{-1}\right)^{-1}+
O\left(q^{-d/2}\sigma(\widetilde{h})\right)\,.\]
Since $f_i$'s are always prime powers, $\sigma(\widetilde{h}) \le 2^r$. The total contribution to $\langle(\triangle_Z)^r \rangle$ from the error term $O\left(q^{-d/2}\sigma(\widetilde{h})\right)$ is bounded by
\[T_2 \le \sum_{n_1, \ldots, n_r \le Z} \prod_{i=1}^r q^{-n_i}n_i^{-1} \sum_{\substack{\deg f_i =n_i \\ 1 \le i \le r}} \Lambda(f_1) \cdots \Lambda(f_r) q^{-d/2}2^r.\]
This can be estimated as
\begin{eqnarray} \label{4:t2} T_2 \le q^{-d/2}2^r (1+\log Z)^r\,.\end{eqnarray}
The total contribution from the main term $\prod_{P|\widetilde{h}}\left(1+|P|^{-1}\right)^{-1}$ is \begin{eqnarray*} \sum_{n_1,\ldots,n_r \le Z} \prod_{i=1}^r q^{-n_i}n_i^{-1} \sum_{\substack{\deg f_i =n_i \\ 1 \le i \le r\\ f_1 \cdots f_r=h^2}} \Lambda(f_1) \cdots \Lambda(f_r)\prod_{P|h}\left(1+|P|^{-1}\right)^{-1}\,.\end{eqnarray*}
Removing the restrictions that $n_1, \ldots, n_r \le Z$, all terms being non-negative, we can bound the above sums as
\begin{eqnarray} \label{hr} H(r)=\sum_{n_1,\ldots,n_r \ge 1} \prod_{i=1}^r q^{-n_i}n_i^{-1} \sum_{\substack{\deg f_i =n_i \\ 1 \le i \le r\\ f_1 \cdots f_r=h^2}} \Lambda(f_1) \cdots \Lambda(f_r)\prod_{P|h}\left(1+|P|^{-1}\right)^{-1}\,.\end{eqnarray}

Proposition \ref{6:lem3} which we will prove in Section 6 provides us with the estimate
\begin{eqnarray*} H(r) \le C \left(\frac{4 r \log \log r}{\sqrt{q} \log r}\right)^r\end{eqnarray*}
for any positive integer $r \ge 4$, where $C>0$ is an absolute constant. Combining it with (\ref{4:t1}) and (\ref{4:t2}) we obtain
\[\langle(\triangle_Z)^r \rangle \le q^{-d/2} 2^{(Z+1)r}+q^{-d/2}2^r (1+\log Z)^r+C \left(\frac{r}{\sqrt{q \log r}}\right)^r\,.\]
Choosing
\begin{eqnarray} \label{4:z} Z=\left[\frac{d}{\left(\log d\right)^2}\right], \end{eqnarray}
we find that
\begin{eqnarray} \label{4:r-tri}\langle(\triangle_Z)^r \rangle \le C \left(\frac{r}{\sqrt{q \log r}}\right)^r\,\end{eqnarray}
for some absolute constant $C$, where the positive integer $r$ is in the range (\ref{4:r}).

\subsection{Proof of Theorem \ref{1:thm2}}

Since $N_F=\triangle_Z(F)+\epsilon_Z(F)$, we have
\[\langle (N_F)^r \rangle=\langle (\triangle_Z)^r \rangle +E_{Z,r},\]
where
\[E_{Z,r}=\sum_{l=0}^{r-1 }\binom{r}{l} \left\langle (\epsilon_Z)^{r-l}(\triangle_Z)^l\right\rangle\,.\]
Using (\ref{4:e}), (\ref{4:r}), (\ref{4:z}) and (\ref{4:r-tri}) we find that
\[E_{Z,r} \le r^{r+1} \left(\frac{10g}{Z} \,\,q^{-Z/2}\right)\left\langle (\triangle_Z)^{2r}\right\rangle^{1/2} \le C \left(\frac{r}{\sqrt{q \log r}}\right)^r\,\]
for some absolute constant $C$. Therefore again we obtain
\begin{eqnarray*} \label{4:nfc} \left\langle (N_F)^r \right\rangle \le C \left(\frac{r}{\sqrt{q \log r}}\right)^r,\end{eqnarray*}
and the above inequality holds true for any positive even integer $r$ in the range $4 \le r \le \log d$ and $d \ge 100$.

For any $\psi$ with $4 \le \psi \le \log d$, denote
\[M_{\psi}=\#\left\{F \in \HH_{d,q}: |N_F| \ge \psi\right\}.\]
Then for any positive even integer $r$ in the range (\ref{4:r}) we have
\[\frac{\# M_{\psi}}{\#\HH_{d,q}} \,\,\psi^r \le \left\langle (N_F)^r \right\rangle \le C \left(\frac{r}{\sqrt{q \log r}}\right)^r\,.\]
Choosing $r \thickapprox \psi$ we find that
\[\frac{\# M_{\psi}}{\#\HH_{d,q}} \le C \left(\frac{1}{\sqrt{q \log \psi}}\right)^{\psi} \ll \exp \left(-\frac{\psi}{2} \log (q \log \psi)\right)\,.\]
This completes the proof of Theorem \ref{1:thm2}. \quad $\square$

\section{Proof of Theorem \ref{1:thm3}}

The proof of Theorem \ref{1:thm3} is similar to that of Theorem \ref{1:thm2}, except that in the proof of Theorem \ref{1:thm3}, we need exact asymptotic formulas for each fixed $r$-th moment under the limit that $g \to \infty$ and both $g,q \to \infty$, instead of upper bounds as in the proof of Theorem \ref{1:thm2}. We summarize part of the proof of Theorem \ref{1:thm2} as follows. From now on the constants implied by the notation ``$\ll$'' and ``$O$'' may depend on the fixed positive integer $r$.

As $d \to \infty$ or $d,q \to \infty$, let $g=\left[\frac{d-1}{2}\right] \to \infty$. Choose
\begin{eqnarray} \label{5:z}
Z=\left[\frac{d}{(\log d)^2}\right]\,.
\end{eqnarray}
Denote
\begin{eqnarray} \label{5:nf} N_F=\log \# \J_{C_F}-g \log q+\delta_{d/2} \log \left(1-q^{-1}\right).\end{eqnarray}
Then
\begin{eqnarray*} N_F=\triangle_Z(F)+\epsilon_Z(F),\end{eqnarray*}
where
\begin{eqnarray*}
\triangle_Z(F)=\sum_{n \le Z}q^{-n} n^{-1} \sum_{\deg f=n} \Lambda(f) \left(\frac{F}{f}\right),
\end{eqnarray*}
and
\begin{eqnarray} \label{5:e}
|\epsilon_Z(F)| \le \frac{10g}{Z} \,\, q^{-Z/2},\quad \left|\triangle_Z(F)\right| \ll \log Z\,.
\end{eqnarray}
For each positive integer $r$,
\[\langle (\triangle_Z)^r\rangle=\sum_{n_1,\ldots,n_r \le Z} \prod_{i=1}^r q^{-n_i}n_i^{-1} \sum_{\substack{\deg f_i =n_i \\ 1 \le i \le r\\ f_1 \cdots f_r=h^2}} \Lambda(f_1) \cdots \Lambda(f_r)\prod_{P|h}\left(1+|P|^{-1}\right)^{-1}+T_1+T_2, \]
where
\[|T_1|+|T_2| \le q^{-d/2}2^{(Z+1)r} + q^{-d/2}2^r \left(1+\log Z\right)^r \ll q^{-d/3}.\]
We can write the main term as
\[\sum_{h} \prod_{P|h}\left(1+|P|^{-1}\right)^{-1} |h|^{-2} \sum_{\substack{\deg f_i \le Z \\ 1 \le i \le r\\ f_1 \cdots f_r=h^2}} \frac{\Lambda(f_1) \cdots \Lambda(f_r)}{(\deg f_1) \cdots (\deg f_r)}\,.\]
Removing the restriction that $\deg f_1, \ldots, \deg f_r \le Z$ results in an error bounded by
\[\sum_{\substack{h\\ \deg h > Z/2}} \prod_{P|h}\left(1+|P|^{-1}\right)^{-1} |h|^{-2} \sum_{\substack{f_1,\ldots,f_r \\ f_1 \cdots f_r=h^2}} \frac{\Lambda(f_1) \cdots \Lambda(f_r)}{(\deg f_1) \cdots (\deg f_r)}\,.\]
Noticing that $\frac{\Lambda(f_i)}{\deg f_i} \le 1$ and $f_i$'s are all prime powers, the sum over $h$ is actually over all monic polynomials $h \in F[X]$ with $\omega(h) \le r$ and $\deg h >Z/2$, where $\omega(h)$ is the function counting the number of distinct prime factors of $h$. If such an $h$ is chosen, the number of choices for each $f_i$ dividing $h$ which is a prime power is less than $2r \deg h$. Hence the error by removing the restriction that $\deg f_1, \ldots, \deg f_r \le Z$ is bounded by
\[T_3 \le \sum_{\deg h >Z/2} |h|^{-2} (2r \deg h)^r =\sum_{n >Z/2} q^{-n}(2rn)^r \ll q^{-Z/4}\,.\]
Combining these estimates together we obtain
\[\langle (\triangle_Z)^r\rangle=H(r)+T, \]
where $H(r)$ is the same function given in (\ref{hr}) as in the proof of Theorem \ref{1:thm2} and $T \ll q^{-Z/4}$. As in the proof of Theorem \ref{1:thm2}, we write
\[\langle (N_F)^r \rangle=\langle (\triangle_Z)^r \rangle +E_{Z,r},\]
where
\[E_{Z,r}=\sum_{l=1}^r \binom{r}{l} \left\langle (\epsilon_Z)^l(\triangle_Z)^{r-l}\right\rangle \ll q^{-Z/4}\,.\]
Using (\ref{5:z}) and (\ref{5:e}) we find that
\begin{eqnarray} \label{5:nfr} \langle (N_F)^r \rangle=H(r)+O\left(q^{-Z/4}\right).\end{eqnarray}

If $q$ is fixed and $d \to \infty$, then for each fixed $r$,
\[\lim_{d \to \infty}\langle (N_F)^r \rangle=H(r). \]
Now suppose that $X$ is a random variable with
\begin{eqnarray} \label{4:x} \E(X^r) =H(r), \quad \forall r \in \N.\end{eqnarray}
For any $t \in \R$, the characteristic function $\phi(t)$ of $X$ is given by
\[\phi(t)=\E\left(e^{itX}\right).\]
Expanding $e^{itX}$ by using the identity
\begin{eqnarray} \label{4:eexpand} e^x=1+\sum_{n=1}^{\infty}\frac{x^n}{n!}\,,\end{eqnarray}
using (\ref{4:x}) and applying Proposition \ref{6:lem0} for $H(r)$ which we will prove in Section 6, we find that
\[\phi(t)=1+\sum_{n=1}^{\infty} \frac{(it)^n}{n!}\sum_{r=1}^{\infty} \frac{n!}{2^r r!} \sum_{\substack{\lambda_1+\cdots+\lambda_r=n\\ \lambda_i \ge 1}} \sum_{\substack{P_1,\ldots,P_r\\ \mbox{\tiny distinct}}} \prod_{j=1}^r \frac{u_{P_j}^{\lambda_j}+(-1)^{\lambda_j}v_{P_j}^{\lambda_j}}{\lambda_j! \left(1+|P_j|^{-1}\right)}\,,\]
where for any $P \in \F[X]$,
\[u_P=-\log\left(1-|P_j|^{-1}\right), \quad v_P=\log\left(1+|P_j|^{-1}\right).\]
Changing the order of summation again we obtain
\[\phi(t)=1+\sum_{r=1}^{\infty} \frac{1}{2^r r!} \sum_{\substack{P_1,\ldots,P_r\\ \mbox{\tiny distinct}}} \prod_{j=1}^r \left(\sum_{\lambda_j=1}^{\infty} \frac{(it)^{\lambda_j} \left(u_{P_j}^{\lambda_j}+(-1)^{\lambda_j}v_{P_j}^{\lambda_j}\right)}{\lambda_j! \left(1+|P_j|^{-1}\right)} \right)\,.\]
The identity (\ref{4:eexpand}) implies that
\[\phi(t)=1+\sum_{r=1}^{\infty} \frac{1}{2^rr!} \sum_{\substack{P_1,\ldots,P_r\\ \mbox{\tiny distinct}}} \prod_{j=1}^r \left( \frac{\left(1-|P_j|^{-1}\right)^{-it}+\left(1+|P_j|^{-1}\right)^{-it}-2}{\lambda_j! \left(1+|P_j|^{-1}\right)} \right)\,.\]
This completes the proof of (1) of Theorem \ref{1:thm3}.

For (2) of Theorem \ref{1:thm3}, assume that both $d,q \to \infty$. Proposition \ref{6:lem4} which we will prove in Section 6 shows that for each fixed positive integer $r$ we have
\[H(r)= \frac{\delta_{r/2} r!}{2^{r/2} (r/2)!} q^{-r/2}+O\left(q^{-(r+1)/2}\right)\,.\]
Here the implied constant in the notation ``$O$'' depends on $r$. From (\ref{5:nfr}) we find that for any fixed positive integer $r$
\begin{eqnarray*} \langle (N_F)^r\rangle=\frac{\delta_{r/2} r!}{2^{r/2} (r/2)!} q^{-r/2}+O\left(q^{-(r+1)/2}+q^{-Z/4}\right). \end{eqnarray*}
We rewrite it as
\[\left\langle \left(\sqrt{q}N_F\right)^r \right\rangle=\frac{\delta_{r/2} r!}{2^{r/2} (r/2)!} +O\left(q^{-1/2}+q^{-(Z-2r)/4}\right).\]
Using the definition of $N_F$ in (\ref{5:nf}), letting $q,d \to \infty$, also noticing that
\[\lim_{q \to \infty}\sqrt{q}\log\left(1-q^{-1}\right)=0,\]
we find that all moments of $\sqrt{q}\left(\log \#J_{C_F}-g \log q\right)$ as $F$ varies in $\HH_d$ are asymptotic to the corresponding moments of a standard Gaussian distribution, where the odd moments vanish and the even moments are
\[\frac{1}{\sqrt{2 \pi}} \int_{\infty}^{\infty} x^{2r}e^{-x^2/2} \ud \, x = \frac{(2r)!}{2^r r!}\,.\]
This implies that as $q,d \to \infty$ and $F$ varies in the set $\HH_d$, the value $\sqrt{q}\left(\log \#J_{C_F}-g \log q\right)$, considered as a random variable on the space $\HH_d$, converges weakly to a standard Gaussian variable. This completes the proof (2) of Theorem \ref{1:thm3}. $\square$

\section{Analysis of $H(s)$}

\subsection{Proposition 1} Let $\F$ be a finite field of cardinality $q$. For any positive integer $s$, denote
\begin{eqnarray*} H(s)=\sum_{n_1,\ldots,n_s \ge 1} \prod_{i=1}^s q^{-n_i}n_i^{-1} \sum_{\substack{\deg f_i =n_i \\ 1 \le i \le s\\ f_1 \cdots f_s=h^2}} \Lambda(f_1) \cdots \Lambda(f_s)\prod_{P|h}\left(1+|P|^{-1}\right)^{-1}\,.\end{eqnarray*}
We first derive another representation of $H(s)$ which is useful in the proofs of Theorems \ref{1:thm2} and \ref{1:thm3}.

\begin{prop} \label{6:lem0} For any positive integer $s \ge 1$ we have
\begin{eqnarray*} H(s)=\sum_{r=1}^s \frac{s!}{2^r r!}\sum_{\substack{\lambda_1+\cdots+\lambda_r=s\\ \lambda_i \ge 1}}\sum_{\substack{P_1,\ldots,P_r\\ \mbox{\tiny distinct}}} \prod_{i=1}^r\,\frac{u_{P_i}^{\lambda_i}+(-1)^{\lambda_i}
v_{P_i}^{\lambda_i}}{\lambda_i! \left(1+|P_i|^{-1}\right)} \,,\end{eqnarray*}
where the sum on the right side is over all positive integers $\lambda_1,\ldots,\lambda_r$ such that $\lambda_1+\cdots+\lambda_r=s$ and over all distinct monic irreducible polynomials $P_1,\ldots,P_r \in \F[X]$, and
\begin{eqnarray} \label{6:upvp} u_P=-\log \left(1-|P|^{-1}\right), \quad v_P=\log \left(1+|P|^{-1}\right), \quad \forall P \in \F[X]. \end{eqnarray}
\end{prop}

\noindent {\bf Proof.} We rewrite $H(s)$ as
\[H(s)=\sum_{h} \prod_{P|h}\left(1+|P|^{-1}\right)^{-1} |h|^{-2} \sum_{\substack{f_1, \ldots,f_s \\ f_1 \cdots f_s=h^2}} \frac{\Lambda(f_1) \cdots \Lambda(f_s)}{(\deg f_1) \cdots (\deg f_s)}\,.\]
Since $f_i$'s are prime powers, the sum over $h$ is actually over all monic polynomials $h \in \F[X]$ with $\omega(h) \le r$, where $\omega(h)$ is the number of distinct prime factors of $h$. Hence
\begin{eqnarray} \label{6:hss} H(s)=\sum_{r=1}^s H(s,r),\end{eqnarray}
where
\[H(s,r)=\sum_{\substack{h\\ \omega(h)=r}} \prod_{P|h}\left(1+|P|^{-1}\right)^{-1} |h|^{-2} \sum_{\substack{f_1, \ldots,f_s \\ f_1 \cdots f_s=h^2}} \frac{\Lambda(f_1) \cdots \Lambda(f_s)}{(\deg f_1) \cdots (\deg f_s)}\,.\]
If $\omega(h)=r$, write explicitly $h=P_1^{a_1}\cdots P_r^{a_r}$ for some distinct primes $P_1,\ldots,P_r$ and exponents $a_1,\ldots,a_r \ge 1$, then
\[H(s,r)=\frac{1}{r!}\sum_{\substack{P_1,\ldots,P_r\\ \mbox{\tiny distinct}}} \sum_{\substack{a_1,\ldots,a_r \ge 1\\ h=P_1^{a_1}\cdots P_r^{a_r}}}\prod_{i=1}^r\left(1+|P_i|^{-1}\right)^{-1}|P_i|^{-2a_i} \sum_{\substack{f_1, \ldots,f_s \\ f_1 \cdots f_s=h^2}} \frac{\Lambda(f_1) \cdots \Lambda(f_s)}{(\deg f_1) \cdots (\deg f_s)}\,.\]
Since each $f_i$ is a prime power and $f_1 \cdots f_s=P_1^{2a_1}\cdots P_r^{2a_r}$, there are finitely many ways to assign prime powers to each $f_i$, according to which we will break $H(s,r)$ into many subsums. With that in mind, for each partition of the set of indexes
\[\{1,2,\ldots,s\}=\bigcup_{i=1}^r A_i,\quad \#A_i=\lambda_i \ge 1, \forall i,\]
it satisfies the property that
\[\sum_{i=1}^r\lambda_i=s\,.\]
We say $(A_1,\ldots,A_r)$ is the type of $(f_1,\ldots,f_r)$ with $f_1 \cdots f_r=h^2$, namely whenever $j \in A_i$, then $f_j$ is a power of $P_i$. Suppose that $f_i=Q_i^{e_i}$ for some prime $Q_i \in \{P_1,\ldots,P_r\}$ and exponent $e_i \ge 1$, and the type of $(f_1,\ldots,f_r)$ is $(A_1,\ldots,A_r)$, since $f_1 \cdots f_s=P_1^{2a_1}\cdots P_r^{2a_r}$, comparing the exponents of $P_j$ on both sides we find that
\begin{eqnarray} \label{6:par} \sum_{i \in A_j}e_i=2a_j \,\, \forall 1 \le j \le r, \end{eqnarray}
and
\[\frac{\Lambda(f_1) \cdots \Lambda(f_s)}{(\deg f_1) \cdots (\deg f_s)}=\frac{1}{e_1 \cdots e_s}.\]
Instead of summing over all integers $a_1,\ldots,a_r$, we sum over all positive integers $e_1,\ldots,e_s$ which satisfy the conditions (\ref{6:par}). Noticing that the value only depends on the vector of integers $(\lambda_1,\ldots,\lambda_r)$ such that
\[\sum_{i=1}^r \lambda_i=s,\]
hence we can write $H(s,r)$ as
\[H(s,r)=\frac{s!}{r!}\sum_{\substack{\lambda_1+\cdots+\lambda_r=s\\ \lambda_i \ge 1}}\sum_{\substack{P_1,\ldots,P_r\\ \mbox{\tiny distinct}}} \prod_{i=1}^r\left(\frac{\left(1+|P_i|^{-1}\right)^{-1}}{\lambda_i !} \sum_{\substack{a_1+ \cdots+a_{\lambda_i} \equiv 0 \pmod{2} \\ a_j \ge 1}} \frac{|P_i|^{-a_1-\cdots-a_{\lambda_i}}}{a_1 \cdots a_{\lambda_i}}\right)\,.\]

For each prime $P$ and positive integer $\lambda$, denote
\[\eta(\lambda)=\eta_P(\lambda):=\sum_{\substack{a_1+ \cdots+a_{\lambda} \equiv 0 \pmod{2} \\ a_i \ge 1}} \frac{|P|^{-a_1-\cdots-a_{\lambda}}}{a_1 \cdots a_{\lambda}}\,,\]
and
\[\tau(\lambda)=\tau_P(\lambda):=\sum_{\substack{a_1+ \cdots+a_{\lambda} \equiv 1 \pmod{2} \\ a_i \ge 1}} \frac{|P|^{-a_1-\cdots-a_{\lambda}}}{a_1 \cdots a_{\lambda}}\,.\]
Since
\[-\log (1-x)=\sum_{n \ge 1} \frac{x^n}{n}, \quad |x|<1, \]
we find
\begin{eqnarray} \label{6:e1}
\eta(1)=- \frac{1}{2}\log \left(1-|P|^{-2}\right),
\end{eqnarray}
and
\begin{eqnarray} \label{6:eta1} \eta(\lambda)+\tau(\lambda)=\sum_{a_1, \cdots a_{\lambda} \ge 1 } \frac{|P|^{-a_1-\cdots-a_{\lambda}}}{a_1 \cdots a_{\lambda}}=(-1)^{\lambda}\log^{\lambda}\left(1-|P|^{-1}\right)\,.
\end{eqnarray}
Combining (\ref{6:e1}) and (\ref{6:eta1}) we have
\begin{eqnarray*}
\tau(1)=-\log \left(1-|P|^{-1}\right)+ \frac{1}{2}\log \left(1-|P|^{-2}\right).
\end{eqnarray*}

For $\lambda \ge 2$, we can write
\[\eta(\lambda)=\sum_{\substack{a_2+ \cdots+a_{\lambda} \equiv 0 \pmod{2} \\ a_i \ge 1}} \left(\prod_{i=1}^{\lambda}\frac{|P|^{-a_i}}{a_i}\right) \eta(1)+\sum_{\substack{a_2+ \cdots+a_{\lambda} \equiv 1 \pmod{2} \\ a_i \ge 1}} \left(\prod_{i=1}^{\lambda}\frac{|P|^{-a_i}}{a_i}\right) \tau(1).\]
This shows that
\begin{eqnarray} \label{6:eta2} \eta(\lambda)=\eta(1) \eta(\lambda-1)+\tau(1) \tau(\lambda-1).\end{eqnarray}
Similarly for $\lambda \ge 2$,
\begin{eqnarray} \label{6:eta3} \tau(\lambda)=\eta(1) \tau(\lambda-1)+\tau(1) \eta(\lambda-1).\end{eqnarray}
We can assign the initial values
\[\eta(0)=1, \quad \tau(0)=0,\]
so that the recursive relations (\ref{6:eta2}) and (\ref{6:eta3}) hold for any $\lambda \ge 1$. Subtracting these two recursive relations we obtain
\[\eta(\lambda)-\tau(\lambda)=\left(\eta(1)-\tau(1)\right)
\left(\eta(\lambda-1)-\tau(\lambda-1)\right).\]
Applying this relation recursively and using (\ref{6:eta1}) we conclude that
\[\eta(\lambda)=\frac{1}{2} \left(u_P^{\lambda}+(-1)^{\lambda}v_P^{\lambda}\right),\]
where
\begin{eqnarray*} u_P=-\log \left(1-|P|^{-1}\right), \quad v_P=\log \left(1+|P|^{-1}\right). \end{eqnarray*}
Therefore $H(s,r)$ can be written as
\begin{eqnarray*} H(s,r)=\frac{s!}{2^r r!}\sum_{\substack{\lambda_1+\cdots+\lambda_r=s\\ \lambda_i \ge 1}}\sum_{\substack{P_1,\ldots,P_r\\ \mbox{\tiny distinct}}} \prod_{i=1}^r\,\frac{u_{P_i}^{\lambda_i}+(-1)^{\lambda_i}
v_{P_i}^{\lambda_i}}{\lambda_i! \left(1+|P_i|^{-1}\right)}\,.\end{eqnarray*}
Returning to (\ref{6:hss}) completes the proof of Proposition \ref{6:lem0}. \quad $\square$

\subsection{Proposition 2} We will prove the following.
\begin{prop} \label{6:lem3} For any positive integer $s \ge 4$ we have
\[H(s) \le C \left(\frac{4s \log \log s}{\sqrt{q}\log s}\right)^s\,,\]
where $C>0$ is an absolute constant.
\end{prop}

\noindent {\bf Proof.} Denote for each positive integer $\lambda$
\[h(\lambda)=\sum_{P}\left(u_{P}^{\lambda}+(-1)^{\lambda}
v_{P}^{\lambda}\right), \]
where the summation is over all monic irreducible polynomials $P \in \F[X]$. An upper bound of $H(s,r)$ is given by
\begin{eqnarray} \label{6:gsr} G(s,r) = \frac{s!}{2^r r!}\sum_{\substack{\lambda_1+\cdots+\lambda_r=s\\ \lambda_i \ge 1}} \prod_{i=1}^r \frac{h(\lambda_i)}{\lambda_i!}\,.\end{eqnarray}

\begin{lemma} \label{6:lem1} For each positive integer $\lambda \ge 2$, we have
\[h(2+\lambda) < h(2)h(\lambda)\,.\]
\end{lemma}
\noindent {\bf Proof.} We have
\[h(2)h(\lambda)=\sum_{P,Q} \left(u_{P}^{2}+
v_{P}^{2}\right)\left(u_{Q}^{\lambda}+(-1)^{\lambda}
v_{Q}^{\lambda}\right),\]
hence
\[h(2)h(\lambda)>\sum_{P=Q} \left(u_{P}^{2}+
v_{P}^{2}\right)\left(u_{P}^{\lambda}+
(-1)^{\lambda}v_P^{\lambda}\right).\]
Since $u_P>v_P>0$ for any $P$, we have
\[\left(u_{P}^{2}+
v_{P}^{2}\right)\left(u_{P}^{\lambda}+
(-1)^{\lambda}v_P^{\lambda}\right)-\left(u_{P}^{2+ \lambda}+(-1)^{\lambda}
v_{P}^{2+\lambda}\right)=u_P^2v_P^2\left(u_P^{\lambda-2}+
(-1)^{\lambda}v_P^{\lambda-2}\right) \ge 0,\]
hence
\[h(2)h(\lambda)>\sum_{P} \left(u_{P}^{2+ \lambda}+(-1)^{\lambda}
v_{P}^{2+\lambda}\right)=h(2+\lambda).\]
This completes the proof of Lemma \ref{6:lem1}. \quad $\square$

\begin{lemma} \label{6:lem2}
\[h(3) < h(2)^{3/2}.\]
\end{lemma}
\noindent {\bf Proof.} Since
\[h(2)^3=\sum_{P,Q,R}\left(u_P^2+v_P^2\right)
\left(u_Q^2+v_Q^2\right)\left(u_R^2+v_R^2\right),\]
we have
\[h(2)^3>\sum_{P,Q=R}\left(u_P^2+v_P^2\right)
\left(u_Q^2+v_Q^2\right)^2. \]
Similarly,
\[h(2)^3>\sum_{P=R,Q}\left(u_P^2+v_P^2\right)^2
\left(u_Q^2+v_Q^2\right). \]
Noting that
\[h(3)^2=\sum_{P,Q} \left(u_P^3-v_P^3\right)\left(u_Q^3-v_Q^3\right)<\sum_{P,Q} u_P^3u_Q^3\,,\]
and
\[u_P^3u_Q^3 \le \frac{1}{2}\left( u_P^2u_Q^4+u_P^4u_Q^2\right)<\frac{1}{2}\left( \left(u_P^2+v_P^2\right)\left(u_Q^2+v_Q^2\right)^2+
\left(u_P^2+v_P^2\right)^2\left(u_Q^2+v_Q^2\right)\right),\]
summing over all $P$ and $Q$ we find that
\[h(3)^2 < h(2)^3.\]
This completes the proof of Lemma \ref{6:lem2}. \quad $\square$

\begin{lemma} \label{6:lem22}
\[h(1)<h(2)^{1/2}.\]
\end{lemma}
\noindent {\bf Proof.} This can be checked explicitly. First
\[h(1)=\sum_{P} \log \frac{1}{1-|P|^{-2}}=\log \prod_{P}\left(1-|P|^{-2}\right)^{-1}\,.\]
From the zeta function of the rational function field $K=\F(X)$ we know that
\[\prod_{P}\left(1-|P|^{-2}\right)^{-1}=\left(1-q^{-1}\right)^{-1},\]
hence
\[h(1)=-\log \left(1-q^{-1}\right)\,.\]

On the other hand,
\begin{eqnarray} \label{6:h2} h(2)=\sum_{P}\log^2\left(\frac{1}{1-|P|^{-1}}\right)+
\log^2\left(1+|P|^{-1}\right).\end{eqnarray}
The terms with $\deg P=1$ (there are $q$ of them) already contribute
\[q \left(\log^2\left(\frac{1}{1-q^{-1}}\right)+
\log^2\left(1+q^{-1}\right)\right)>q\log^2\left(\frac{1}{1-q^{-1}}\right)\]
to $h(2)$. This completes the proof of Lemma \ref{6:lem22}. \quad $\square$

From Lemma \ref{6:lem2} and Lemma \ref{6:lem22} we know that
\[h(1)^2h(3)^2<h(2)h(2)^3=h(2)^4,\]
hence
\begin{eqnarray} \label{6:h13} h(1)h(3) <h(2)^2.\end{eqnarray}

Suppose that $s$ is a positive integer. For any positive integers $\lambda_1,\ldots,\lambda_r$ such that
\[\sum_{i=1}^r \lambda_i=s,\]
from Lemma \ref{6:lem1}--\ref{6:lem22} and using (\ref{6:h13}), we find that
\[\prod_{i=1}^r h(\lambda_i) \le h(2)^{s/2}\,.\]
Using this in (\ref{6:gsr}) to get an upper bound for $H(r,s)$ and then returning to $H(s)$ in (\ref{6:hss}), we obtain that for any positive integer $s$,
\[H(s) \le \sum_{r=1}^s \frac{s!h(2)^{s/2}}{2^r r!}\sum_{\substack{\lambda_1+\cdots+\lambda_r=s\\ \lambda_i \ge 1}} \frac{1}{\lambda_1! \cdots \lambda_r!}\,.\]

From definition of $h(2)$ in (\ref{6:h2}) and using
\[\log (1+x) \le x,\quad 0<x<1,\]
we find
\[h(2) \le \sum_{P} \frac{1}{\left(|P|-1\right)^2}+\frac{1}{|P|^2}\,. \]
Summing over all monic polynomials $h \in \F[X]$ instead of monic irreducible polynomials $P \in \F[X]$, we obtain
\[h(2)< \sum_{n=1}q^n \left\{\frac{1}{\left(q^n-1\right)^2}+\frac{1}{q^{2n}}\right\}\,,\]
hence
\[h(2)<\left(\left(1-q^{-1}\right)^{-2}+1\right)\sum_{n \ge 1}q^{-n} \le 10 q^{-1}\,.\]
Also from the identity
\[\left(x_1+\cdots+x_r\right)^s=
\sum_{\substack{\lambda_1+\cdots+\lambda_r=s\\ \lambda_i \ge 0}}
\frac{s!}{\lambda_1! \cdots \lambda_r!}\,\,\,x_1^{\lambda_1} \cdots x_r^{\lambda_r},\]
we find that
\[\sum_{\substack{\lambda_1+\cdots+\lambda_r=s\\ \lambda_i \ge 1}}
\frac{s!}{\lambda_1! \cdots \lambda_r!} < r^s. \]
Therefore
\[H(s)<\sum_{r=1}^s \frac{1}{2^r r!}\, \left(\frac{\sqrt{10}\,r}{\sqrt{q}}\right)^s\,.\]
To find an upper bound, denote
\[a_r=\frac{r^s}{2^r r!}.\]
Then
\[\frac{a_{r+1}}{a_r}=\frac{1}{2r} \left(1+\frac{1}{r}\right)^{s-1}\,.\]
If $s \ge 100$, we choose
\[l_0=\left[\frac{s \log \log s}{\log s}\right].\]
We find that
\[\sum_{r=1}^{l_0}\frac{1}{2^r r!}\, \left(\frac{\sqrt{10}\,r}{\sqrt{q}}\right)^s \le \left(\frac{\sqrt{10}\,s \log \log s}{\sqrt{q}\log s}\right)^s \sum_{r=1}^{\infty}\frac{1}{2^r r!} =e^{1/2}\left(\frac{4s \log \log s}{\sqrt{q}\log s}\right)^s\,. \]
For the choice of $l_0$ we have
\[l_0 \log l_0 =s \log \log s \left(1+o_s(1)\right) \ge s,\]
and from it we derive
\[\frac{a_{r+1}}{a_r} \le \frac{1}{2}, \quad \forall r \ge l_0. \]
Hence
\[\sum_{r=l_0}^{s}\frac{1}{2^r r!}\, \left(\frac{\sqrt{10}\,r}{\sqrt{q}}\right)^s \le \frac{1}{2^{l_0} {l_0}!}\, \left(\frac{\sqrt{10}\,l_0}{\sqrt{q}}\right)^s\left(1+2^{-1}+2^{-2}+\cdots\right) \le \left(\frac{\sqrt{10}\,s \log \log s}{\sqrt{q}\log s}\right)^s\,. \]
Therefore
\[H(s) \le 3 \left(\frac{\sqrt{10}\,s \log \log s}{\sqrt{q}\log s}\right)^s\]
for any positive integer $s \ge 100$. On the other hand, if $4 \le s <100$, we use the trivial estimate
\[H(s)<\left(\frac{\sqrt{10}\,s}{\sqrt{q}}\right)^s
\sum_{r=0}^{\infty}\frac{1}{2^r r!}=e^{1/2}\left(\frac{\sqrt{10}\,s}{\sqrt{q}}\right)^s,\]
this again is bounded by $C \left(\frac{\sqrt{10}\,s \log \log s}{\sqrt{q}\log s}\right)^s$ for some absolute constant $C>0$ for $4 \le s <100$. This completes the proof of Proposition \ref{6:lem3}. \quad $\square$

\subsection{Proposition 3} Finally, we prove the following.
\begin{prop} \label{6:lem4} If $s$ is a fixed positive integer, then \[H(s) = \frac{\delta_{s/2}\, s!}{2^{r/2} (s/2)!} \,\, q^{-s/2}+O_s\left(q^{-(s+1)/2}\right)\,,\]
as $q \to \infty$, where for any $\gamma \in \R$,
\[\delta_{\gamma}=\left\{\begin{array}{cc}1&\gamma \in \Z,\\
0& \gamma \not \in \Z.\end{array}\right.\]
\end{prop}
\noindent {\bf Proof.} For any $\lambda \in \N$ and $P \in \F[X]$, denote
\[g(\lambda,|P|)= \frac{u_P^{\lambda}+(-1)^{\lambda}v_P^{\lambda}}{1+|P|^{-1}}, \]
where $u_P$ and $v_P$ are given in (\ref{6:upvp}). We use Taylor series expansions of $-\log (1-x)$ and $\log (1+x)$ ($|x| \le 1/2$) given by
\[-\log (1-x)=x+\frac{x^2}{x}+\frac{x^3}{3}+\cdots=x\left(1+\frac{x}{2}+O(x^2)\right),\]
and
\[\log (1+x)=x-\frac{x^2}{x}+\frac{x^3}{3}+\cdots=x\left(1-\frac{x}{2}+O(x^2)\right)\]
to deduce
\begin{eqnarray} \label{7:g} g(\lambda,|P|)=\left\{\begin{array}{ccc} 2|P|^{-\lambda} \left(1+O\left(|P|^{-1}\right)\right)&:&\lambda \equiv 0 \pmod{2},\\
\lambda|P|^{-\lambda-1} \left(1+O\left(|P|^{-1}\right)\right)&:&\lambda \equiv 1 \pmod{2}.\end{array}\right.\end{eqnarray}
For any $n \in \N$, denote
\[\pi_q(n)=\#\{P \in \F[X]: P \mbox{ is monic, irreducible and } \deg P=n\}\,.\]
It is known (\cite{ros}) that
\begin{eqnarray} \label{7:pi} \pi_q(n)=\frac{q^n}{n} \left(1+O(q^{-n/2})\right). \end{eqnarray}
Fix a positive integer $s$. For any positive integer $\lambda \le s$, denote
\[z(\lambda)=\sum_{P}g(\lambda,|P|)=\sum_{n \ge 1}g(\lambda,q^n) \pi_q(n).\]
If $\lambda$ is even, from (\ref{7:g})
\[z(\lambda)=\sum_{n \ge 1}\frac{2q^{-n \lambda} \left(1+O(q^{-n})\right)}{1+q^{-n}} \pi_q(n)\,.\]
For $n=1$, since $\pi_q(1)=q$, this gives us the value $2q^{-\lambda+1} \left(1+O(q^{-1})\right)$. For $n \ge 2$, using (\ref{7:pi}) and noting that $\lambda \ge 2$, we find that all such terms added together contribute to a value bounded by $O(q^{-\lambda})$. Hence
\begin{eqnarray} \label{7:ze} z(\lambda)=2q^{-\lambda+1} \left(1+O_s(q^{-1})\right), \quad \mbox{if } \lambda \equiv 0 \pmod{2}\,.\end{eqnarray}
The notation ``$O_s$'' in the above formula is to stress that the implied constant may depend on the value $s$ since $1 \le \lambda \le s$. Similarly if $\lambda$ is odd, we obtain that
\begin{eqnarray} \label{7:zo} z(\lambda)=\lambda q^{-\lambda}\left(1+O_s(q^{-1})\right), \quad \mbox{if } \lambda \equiv 1 \pmod{2}\,.\end{eqnarray}

We write
\begin{eqnarray} \label{6:hs} H(s)=\sum_{r=1}^s\frac{s!}{2^r r!}\,\, \sum_{\substack{\lambda_1+\cdots+\lambda_r=s\\ \lambda_i \ge 1}} f(\lambda_1,\ldots,\lambda_r)\,, \end{eqnarray}
where \begin{eqnarray*} f(\lambda_1,\ldots,\lambda_r)= \sum_{\substack{P_1,\ldots,P_r\\ \mbox{\tiny distinct}}} \prod_{i=1}^r\,\,\frac{g(\lambda_i,|P_i|)}{\lambda_i!} \,.\end{eqnarray*}
For any positive integers $\lambda_1,\ldots,\lambda_r$ with
\[\sum_{i=1}^r \lambda_i=s,\]
denote
\[\{1,2,\ldots,r\}=A \bigcup B,\]
where
\[A=\{1 \le i \le r: \lambda_i \mbox{ is even}\},\quad B=\{1 \le i \le r: \lambda_i \mbox{ is odd}\}.\]
Since for any $i \in A$, $\lambda_i \ge 2$, and
\[s=\sum_{i \in A} \lambda_i+\sum_{i \in B} \lambda_i,\]
we check that if $\#B \ge 1$ then $\#A \le (s-1)/2$, or if there is an $i \in A$ with $\lambda_i \ge 4$, then $\#A \le (s-2)/2$. In either of these cases, since
\[f(\lambda_1,\ldots,\lambda_r) \le \left(\prod_{i \in A} z(\lambda_i)\right) \left(\prod_{i \in B} z(\lambda_i)\right), \]
using (\ref{7:ze}) and (\ref{7:zo}) we obtain
\[f(\lambda_1,\ldots,\lambda_r) \ll_s \left(\prod_{i \in A}q^{-\lambda_i+1}\right) \left(\prod_{i \in B} q^{-\lambda_i}\right)=q^{-s+\#A} \le q^{-(s+1)/2}\,.\]

The case that is not covered by the above consideration is that $\#B=0$ and there is no $i \in A$ with $\lambda_i \ge 4$, that is, $\lambda_i=2$ for any $1 \le i \le r$. That means that $s=2r$ is even. The only term left is
\begin{eqnarray} \label{7:m1} f(2,\ldots,2)=\frac{1}{2^r} \sum_{\substack{P_1,\ldots,P_r\\ \mbox{\tiny distinct}} } \prod_{i=1}^r g\left(2,|P_i|\right)\,.\end{eqnarray}

Since
\[\sum_{\substack{P_1,\ldots,P_r\\ \mbox{\tiny not distinct}} } \prod_{i=1}^r g\left(2,|P_i|\right) \ll_s \sum_{\substack{P_1,\ldots,P_r\\ P_1=P_2} } \prod_{i=1}^r g\left(2,|P_i|\right)=\left(\sum_{P}g\left(2,|P|\right)^2\right)
\left(\sum_{P}g\left(2,|P|\right)\right)^{r-2},\]
and we can check easily that
\[\sum_{P}g\left(2,P\right)^2 \ll q^{-3}.\]
Using the above and (\ref{7:ze}) we find that removing the restrictions that $P_1,\ldots,P_r$ are distinct in (\ref{7:m1}) would result in an error bounded by $O_s\left(q^{-r-1}\right)=O_s\left(q^{-(s+1)/2}\right)$. The main term in (\ref{7:m1}) is given by
\[\frac{1}{2^r}\left(\sum_{P}g\left(2,|P|\right)\right)^{r}
=q^{-r}\left(1+O_s(q^{-1})\right)=q^{-s/2}\left(1+O_s(q^{-1})\right)\,.\]
Combining all the above computation together and returning to (\ref{6:hs}), we conclude that
\[H(s) = \frac{\delta_{s/2}\, s!}{2^{r/2} (s/2)!} \,\, q^{-s/2}+O_s\left(q^{-(s+1)/2}\right)\,,\]
as $q \to \infty$. This completes the proof of Proposition \ref{6:lem4}. \quad $\square$


\begin{thebibliography}{99}

\bibitem{ach1} J.D. Achter, \emph{The distribution of class groups of function fields}, J. Pure Appl. Algebra {\bf 204} (2006),  no. 2, 316--333.

\bibitem{ach2} J.D. Achter, \emph{Results of Cohen-Lenstra type for quadratic function fields}, Computational arithmetic geometry, 1--7, Contemp. Math., {\bf 463}, Amer. Math. Soc., Providence, RI, 2008.


\bibitem{adl} L.M. Adleman, M.A. Huang, ``Primality testing and abelian varieties over finite fields'', Lecture Notes in Mathematics, 1512.



\bibitem{buc1} A. Bucur, C. David, B. Feigon, M. Lal\'in, \emph{Statistics for traces of cyclic trigonal curves over finite fields}, International Mathematics Research Notices (2010), 932--967.

\bibitem{buc2} A. Bucur, C. David, B. Feigon, M. Lal\'in, \emph{Biased statistics for traces of cyclic $p$-fold covers over finite fields}, to appear in Proceedings of Women in Numbers, Fields Institute Communications.


\bibitem{deu} M. Deuring. \emph{Die Typen der Multiplikatorenringe elliptischer Funktionenk\"o per}, Abh. Math. Sem. Hansischen Univ., {\bf 14} (1941), 197--272.


\bibitem{dia} P. Diaconis, S. Evans, \emph{Linear functionals of eigenvalues of random matrices}, Trans. Amer. Math. Soc. {\bf 353} (2001), no. 7, 2615--2633.

\bibitem{fai} D. Faifman, Z. Rudnick, \emph{Statistics of the zeros of zeta functions in families of hyperelliptic curves over a finite field}, arXiv:0803.3534. To appear in Compositio Math.



\bibitem{kat} N. M. Katz, P. Sarnak, ``Random Matrices, Frobenius Eigenvalues, and Monodromy'', Amer. Math. Soc. Colloq. Publ., vol. 45, American Mathematical Socitey, Providence, RI, 1999.


\bibitem{kur} P. Kurlberg, Z. Rudnick, \emph{The fluctuations in the number of points on a hyperelliptic curve over a finite field}, J. Number Theory Vol. 129 {\bf 3} (2009), 580-587.

\bibitem{len1} H. W. Lenstra, Jr., J. Pila, C. Pomerance, \emph{A hyperelliptic smoothness test. I}, Philos. Trans. Roy. Soc. London Ser. A {\bf 345} (1993), no. 1676, 397--408.

\bibitem{len2} H. W. Lenstra, Jr., \emph{Factoring integers with elliptic curves}, Ann. of Math. (2) {\bf 126} (1987), no. 3, 649--673.


\bibitem{lor} D. Lorenzini, ``An invitation to arithmetic geometry'', Graduate Studies in Mathematics, 9. American Mathematical Society, Providence, RI, 1996.


\bibitem{mor} C. Moreno, ``Algebraic curves over finite fields'', Cambridge Tracts in Mathematics {\bf 97}, Cambridge University Press, 1991.

\bibitem{que} H.G. Quebbemann. \emph{Estimates of regulators and class numbers in function fields}, J. Reine Angew. Math., {\bf 419} (1991), 79--87.

\bibitem{rosen} M.Y. Rosenbloom, M.A. Tsfasman, \emph{Multiplicative lattices in global fields}, Invent. Math., {\bf 101} (1990), 687--696.

\bibitem{ros} M. Rosen, ``Number theory in function fields''. Graduate Texts in Mathematics, {\bf 210}. Springer-Verlag, New York, 2002.

\bibitem{rud} Z. Rudnick, \emph{Traces of high powers of the Frobenius class in the hyperelliptic ensemble}, Acta Arith. {\bf 143} (2010), 81-99.







\bibitem{shp} I. Shparlinski, \emph{On the size of the Jacobians of curves over finite fields}, Bull. Braz. Math. Soc. (N.S.) {\bf 39} (2008), no. 4, 587--595.


\bibitem{ste} A. Stein, E. Teske. \emph{Explicit bounds and heuristics on class numbers in hyperelliptic function fields}, Math. Comp., {\bf 71} (2002), 837--861.

\bibitem{tsf} M. Tsfasman, \emph{Some remarks on the asymptotic number of points}, Coding theory and algebraic geometry (Luminy, 1991). Lect. Notes in Math., vol. 1518, Springer, (1992), 178--192.

\bibitem{ven} A. Venkatesh, S. Ellenberg, \emph{Statistics of number fields and function fields}, Proceedings of the International Congress of Mathematicians Hyderabad, India, 2010.

\bibitem{wei} A. Weil, \emph{Sur les Courbes Alg\'ebriques et les Vari\'et\'es qui s'en D\'eduisent}, Publ. Inst. Math. Univ. Strasbourg {\bf 7} (1945),  Hermann et Cie., Paris, 1948.


\end{thebibliography}
\end{document}